\documentclass{amsart}
\numberwithin{equation}{section}

\newtheorem{thm}{Theorem}[section]
\newtheorem{lem}[thm]{Lemma}

\newtheorem{defin}[thm]{Definition}

\begin{document}

\author{Ravshan Ashurov}
\address{Ashurov R: Institute of Mathematics of Uzbekistan,
Tashkent, Uzbekistan}
\email{ashurovr@gmail.com}

\author{Yusuf Fayziev}
\address{Fayziev Yu: National University of Uzbekistan,  Tashkent, Uzbekistan}
\email{fayziev.yusuf@mail.ru}

\author{Muattar Khudoykulova}
\address{Khudoykulova M: National University of Uzbekistan,  Tashkent, Uzbekistan}
\email{muattarxudoyqulova2000@gmail.com}

\small

\title[On the non-local problem ...] {On the non-local problem for Boussinesq type fractional equation}

\begin{abstract}

In recent years, the Boussinesq type fractional partial differential equation has attracted much attentions of researchers
for its practical importance. In this paper we study a non-local problem for the Boussinesq type equation $D_t^\alpha u(t)+A D_t^\alpha u(t)+\nu^2A u(t)=0,\,\,  0< t< T,\,\, 1<\alpha<2,$ where $D_t^\alpha$ is the Caputo fractional derivative and $A$ is abstract operator. In the classical case, i.e. at $\alpha=2$, this problem was studied earlier and an interesting effect was discovered: the well-posedness of the problem significantly depends on the length of the time interval and the parameter $\nu$. This note shows that for the case of a fractional equation there is no such effect: the problem is well-posed for any $T$ and $\nu$.

\end{abstract}

\maketitle

\section{Introduction}

Let $H$ be a separabel Hilbert space and operator $A:D(A)\subset H\rightarrow H$ be an arbitrary unbounded, positive self-adjoint operator and we assume that $A$ has a compact inverse operator $A^{-1}$, where $D(A)$ is domain of operator $A$ given below. Let $\lambda_k$ and $\{v_k\}$ be the eigenvalues and corresponding eigenfunctions of operator $A$, that is
\begin{equation}\label{Eigen}
    Av_k = \lambda_kv_k, \,\,\, k=1, 2, \dots \, \, .
\end{equation}
 
Let us introduce  the Caputo fractional derivative $D^\alpha_t$ of order  $\alpha\in(1,2)$  of the vector-valued function $h(t)\in H$ (see, for example \cite{Lizama})

\begin{equation}\label{def0}
D_t^\alpha h(t) = \frac{1}{\Gamma
(2-\alpha)}\int\limits_0^t\frac{h''(\xi)}{(t-\xi)^{\alpha-1}} d\xi, 
\quad t>0,
\end{equation}
provided the right-hand side exists. Here $\Gamma(\alpha)$ is
Euler's gamma function.

Let $1<\alpha<2$. Consider the following fractional differential equation
\begin{equation}\label{prob1}
D_t^\alpha u(t)+AD_t^\alpha u(t)+\nu^2Au(t)=0,\quad   0< t< T,
\end{equation}
with  non-local conditions
\begin{equation}\label{ConB}
u(0)=u(T),
\end{equation}
and
\begin{equation}\label{ConIn}
\int\limits_0^T u(t)dt=\varphi,
\end{equation}
 where $\varphi \in H$ a given vector and $\nu>0$ fixed number.

 Note that since the abstract operator $A$ is required only to have a complete orthonormal system of eigenfunctions, then any elliptic operator can be considered as $A$. For example, if we take $L_2(\Omega )$, $ \Omega\subset \mathbb{R}^N$, as the Hilbert space $H$, then as $A$ we can take the Laplace operator $(-\Delta)$ with the Dirichlet condition.

 The solution to problem (\ref{prob1}) - (\ref{ConIn})  will be understood in the sense of the following definition:
  \begin{defin}\label{def} \,If a function $u(t) \in C([0,T]; H)$, $D_t^\alpha u(t)$, $Au(t)$, $AD_t^\alpha u(t) \in C((0,T) ; H)$ and satisfies all the conditions of problem (\ref{prob1}) - (\ref{ConIn}), then it is called the solution of problem (\ref{prob1}) - (\ref{ConIn}).
\end{defin}

 Equation (\ref{prob1}) has different names for different values of  parameter $\alpha$. So if $\alpha=1$ it is called a Barenblatt-Zheltov-Kochina type differential equation (see \cite{BZhK}) and if $\alpha=2$, it has the name
a Boussinesq type differential equation (see \cite{Bous}). If $0<\alpha<1$, it is called Barenblatt-Zheltov-Kochina type fractional differential equation, in case $1<\alpha<2$, it is called a Boussinesq type fractional differential equation. The Boussinesq type differential equations were introduced by Joseph Boussinesq in 1872 (see \cite{Bous}, eq. 26). Boussinesq's equations are widely used in numerical modeling in coastal engineering to model waves in shallow waters and harbors. Boussinesq's equations are applicable for sufficiently long waves. Although the wave modeling in such cases is perfectly described by the Navier-Stokes equations, it is currently extremely difficult to solve the three-dimensional equations in complex models. Therefore, approximate models such as the Boussinesq's equations can be used to reduce three-dimensional problems to two-dimensional states. (see, for example \cite{Dingeman}-\cite{Kirby}).

There are a number of works  (see, for example, \cite{BZhK}, \cite{BZh} - \cite{AFKh2}), in which specialists consider various initial-boundary value problems for differential and fractional differential equations of the Barenblatt-Zheltov-Kochina type. Since our study relates to the Boussinesq type differential equation, we present some results related specifically to these equations.

Due to the mathematical and physical importance, over the last couple of decades, existence and nonexistence
of solutions of Boussinesq type equations have been extensively studied by many mathematicians and physicists (see, for example \cite{liufangcheng}-\cite{Bjor} and \cite{Jun-Yang}-\cite{Abulwafa} with fractional order, and literature therein). Nonlinear Boussinesq type equations arises in a number of mathematical models of physical processes, for example in the modeling
of surface waves in shallow waters or considering the possibility of energy exchange through lateral surfaces
of the wave guide in the physical study of nonlinear wave propagation in wave guide (see, for example, \cite{ZhangHu} and \cite{Botmart} with fractional order, and literature therein). In \cite{ZhangHu} the authors consider the Cauchy problem of two-dimensional generalized
Boussinesq-type equation $u_{tt}-\Delta u-\Delta u_{tt}+\Delta^2 u+ \Delta f(u) = 0$. Under the assumption
that $f(u)$ is a function with exponential growth at infinity and under some
assumptions on the initial data, the authors prove the existence and in some cases nonexistence of global
weak solution.

Model equations of Boussinesq type (the problem (\ref{prob1}) - (\ref{ConIn}) with $\alpha=2$, $\nu =1$ and $A=-\frac{\partial^2}{\partial x^2}-\frac{\partial^2}{\partial y^2}$, $x,y \in (0, l)$) and equations of mixed type and nonlinear equations containing equations of Boussinesq type are systematically studied in a series of works \cite{AshR14}-\cite{AshR16}. In these works, the existence and uniqueness of the classical solution of initial-boundary value problems was proved and some inverse problems were studied. And in the work \cite{AshR17} problems for the Boussinesq equation with a spectral parameter were studied.

Let us cite two more works \cite{AlKh} and \cite{KhKuch} that motivated the appearance of our research. In these works, the above non-local problem (\ref{prob1}) - (\ref{ConIn}) was studied for a classical partial differential equation in which $A$ is the Laplace operator with the Dirichlet condition. So in the fundamental work \cite{AlKh}, Alimov and Khalmukhamedov studied the following a non-local problem in the cylinder $\Omega\times(0,T)$:
\begin{equation}\label{Aleq}
\left\{
\begin{aligned}
& u_{tt}-\Delta u_{tt}-\nu^2\Delta u=0,\quad x\in \Omega, \quad 0< t< T,\\
&u(x,0)=u(x,T),\quad x\in \Omega, \\
&\int\limits_0^T u(x,t)dt=\varphi(x),
\end{aligned}
\right.
\end{equation}
here $\varphi(x)$ is a given function. The authors discovered an interesting effect: it turns out that the well-posedness of this problem significantly depends on the length of the time interval and on the parameter $\nu$. So if $\frac{\nu T}{2\pi}\in (0,1)$, then the solution exists and is unique for all $\varphi\in D(A)$. Case $\frac{\nu T}{2\pi}\geq 1$ is more complected: if $\frac{\nu T}{2\pi}>1$, and this number is not an natural, then for the existence of a solution it is necessary to require that function $\varphi$ be orthogonal to some eigenfunctions of the Laplace operator, and in this case the solution is not unique. And if the number $\frac{\nu T}{2\pi}$ is an natural, then only orthogonality is not enough; it is necessary that the function $\varphi$ be smoother: $\varphi\in D(A^2)$.

Since the parameter $\nu$ in the equation is fixed, this result means that if the process under study lasts $"$not so long$"$, then a solution to the problem exists for any measurements $\varphi$, and if the process lasts $"$a little$"$  longer, then the solution does not exist for all data $\varphi$.

In the recent work \cite{KhKuch}, problem (\ref{Aleq}) was studied with the kernel $tu(x,t)$  in the integral condition. Similar to the paper \cite{AlKh}, conditions have been found for the time interval $(0,T]$, function $\varphi$ and parameter $\nu$, which guarantees the existence of a solution to the problem.

The question naturally arises: is the effect found in work \cite{AlKh} preserved if instead of the second time-derivative we take the fractional order derivative $\alpha\in (1,2)$ into equation (\ref{prob1})? As shown in the following theorem, the above parameter does not affect the correctness of the problem at all and the solution exists and it is unique for any function $\varphi\in D(A)$.

\begin{thm}\label{w2}
    Let $\varphi\in D(A)$. Then there is a unique solution of problem (\ref{prob1})-(\ref{ConIn}) and it has the form:
    $$
        u(t)=\sum\limits_{k=1}^\infty \left(\frac{\varphi_kE_{\alpha,2}(-\nu_k^2T^\alpha)E_{\alpha,1}(-\nu_k^2t^\alpha)}{T((E_{\alpha,2}(-\nu^2_kT^\alpha))^2+E_{\alpha,3}(-\nu_k^2T^\alpha)(1-E_{\alpha,1}(-\nu_k^2T^\alpha)))}  \right.
        $$
        \begin{equation}\label{yechim}
        +\left. \frac{\varphi_kt(1-E_{\alpha,1}(-\nu_k^2T^\alpha))E_{\alpha,2}(-\nu_k^2t^\alpha)}{T^2((E_{\alpha,2}(-\nu^2_kT^\alpha))^2+E_{\alpha,3}(-\nu_k^2T^\alpha)(1-E_{\alpha,1}(-\nu_k^2T^\alpha)))}\right)v_k,
    \end{equation}
where $\nu_k=\nu\sqrt{\frac{\lambda_k}{1+\lambda_k}}$ and  $\varphi_k=(\varphi,v_k)$ are the Fourier coefficients of function $\varphi$.
\end{thm}

\section{Preliminaries}

In this section, we provide some information about operator $A$ and present new bounds for the Mittag-Leffler function in the case $1<\rho <2$, based on the findings of the study conducted by \cite{Pskhu}.  

The action of the abstract operator $ A $ under consideration on the element $h\in H$ can be written as
$$
A h= \sum\limits_{k=1}^\infty \lambda_k h_k v_k,
$$
where $h_k$ is the Fourier coefficient of the element $h$: $h_k=(h,v_k)$.
 Obviously, the domain of this operator has the
form
$$
D(A)=\{h\in H:  \sum\limits_{k=1}^\infty \lambda_k^{2}
|h_k|^2 < \infty\}.
$$
For elements $h$ and $g$ of $D(A)$ we introduce the norm and inner product as
\[
||h||_1^2=\sum\limits_{k=1}^\infty \lambda_k^{2} |h_k|^2 =
||A h||^2,
\]
\[
(h,g)_1=\sum\limits_{k=1}^\infty \lambda_k^{2} h_k \overline{g}_k,
\]
respectively. Together with this norm $D(A)$ turns into a Hilbert
space.

Recall, the Mittag-Leffler function $E_{\rho,\mu}(t)$ has the form (see e.g. \cite{Gor}, p. 56):
$$
E_{\rho, \mu}(t)= \sum\limits_{n=0}^\infty \frac{t^n}{\Gamma(\rho
n+\mu)},
$$
where $\rho>0$ and $\mu$ complex number.

Next, we establish some two-sided estimates for the Mittag-Leffler function $E_{\rho, \mu}(-t)$, $1<\rho<2$, $t\geq 0$, $\mu= 1, 2, 3, \rho$. The following simple method for obtaining these estimates was suggested to the authors by Professor A.V. Pskhu (see, \cite{Pskhu}).

Let $\phi(\delta, \beta; z)$ stand for the Wright function, defined as
\[
\phi(\delta, \beta; z)= \sum\limits_{k=0}^\infty \frac{z^k}{k! \Gamma(\beta+\delta k)},\,\, \delta>-1, \,\, \beta \in \mathbb{R}, \,\, z\in \mathbb{C}.
\]
Let $1< \alpha <2$. In the work of A.V. Pskhu \cite{Pskhu} for functions $h(t)$ defined at $t\geq 0$, the following integral transform is introduced and studied:
\[
P^{\alpha, \beta} h(t)=t^{\beta-1} \int\limits_0^\infty h(s) \phi\bigg(-\alpha, \beta; -\frac{s}{t^\alpha}\bigg) ds. 
\]
Note that $P^{\alpha, \beta} h(t)$ is some modification of the integral transform introduced by B. Stankovich in 1955 (see \cite{Stankovich}).

Let us present the following statement from \cite{Pskhu}.
\begin{lem}\label{equal} Let $\gamma>0$. Then
\[
P^{\alpha, \beta} t^{\gamma-1}
= t^{\alpha\gamma+\beta-1}\frac{\Gamma(\gamma)}{\Gamma(\alpha \gamma +\beta)}.
\]
\end{lem}

From Lemma \ref{equal}, by the definition of the Mittag-Leffler function we get 
\begin{equation}\label{opert}
    P^{\xi,\eta}[t^{\mu-1}E_{\rho,\mu}(\lambda t^\rho)]=t^{\mu\xi+\eta-1}E_{\rho\xi,\mu\xi+\eta}(\lambda t^{\rho\xi}).
\end{equation}

\begin{lem}\label{theor}(see \cite{Pop}, p. 143.) There is a function $f(\alpha)$ decreasing on the interval $(1, 2)$ such that for any $\alpha \in (1, 2)$ and $\beta > f(\alpha)$ function $E_{\alpha,\beta}(z)$ does not vanish, where $f(\alpha)$ satisfies the following inequality:
$$
\alpha<f(\alpha)<\frac{3}{2}\alpha,\quad 1<\alpha<2.
$$
\end{lem}

\begin{lem}\label{lem1}
Let $\alpha\in(1,2)$. Then the following estimate holds:
$$
|E_{\alpha, 1}(-t^\alpha)|< 1, \quad t>0.
$$
\end{lem}

\emph{Proof.} Let $\mu=1, \rho=1,  \xi=\alpha, \eta=1-\alpha$ and $\lambda=-1$ in equality (\ref{opert}). Then, we have:
$$
E_{\alpha,1}(- t^\alpha)=P^{\alpha,1-\alpha}(E_{1,1}(-t))=P^{\alpha,1-\alpha}e^{-t}.
$$
Using the  inequality $e^{-t}< 1$ and  Lemma \ref{equal} we get:
$$
|E_{\alpha, 1}(-t^\alpha)|< 1, \quad t>0.
$$
Lemma \ref{lem1} is proved.

\begin{lem}\label{lem2}
Let $\alpha\in(1,2)$. Then the following estimate holds:
$$
|E_{\alpha, 2}(-t^\alpha)|<1, \quad t> 0.
$$
\end{lem}
\emph{Proof.} Let $\mu=1, \rho=1,  \xi=\alpha, \eta=2-\alpha$ and $\lambda=-1$. Then from  (\ref{opert}) we have the following equality:
$$
tE_{\alpha,2}(-t^\alpha)=P^{\alpha,2-\alpha}(E_{1,1}(-t))=P^{\alpha,2-\alpha}e^{-t}.
$$
Using the  inequality $e^{-t}< 1$ and Lemma \ref{equal} we get:
$$
|tE_{\alpha, 2}(-t^\alpha)| < P^{\alpha,2-\alpha} 1=t\frac{\Gamma(1)}{\Gamma(2)}=t.
$$
Therefore, for $t>0$, we have
$$
|E_{\alpha, 2}(-t^\alpha)|<1.
$$
 Lemma \ref{lem2} is proved.

\begin{lem}\label{lem3}
Let $\alpha\in(1,2)$. Then the following estimate holds:
$$
0<E_{\alpha, 3}(-t^\alpha)<\frac{1}{2}, \quad t>0.
$$
\end{lem}
\emph{Proof.} Let $\mu=1, \rho=1,  \xi=\alpha, \eta=3-\alpha$ and $\lambda=-1$. Then from  (\ref{opert}) we have the following equality:
$$
E_{\alpha,3}(-t^\alpha)=P^{\alpha,3-\alpha}(E_{1,1}(-t))=P^{\alpha,3-\alpha}e^{-t}.
$$
Using the  inequality $e^{-t}< 1$ and Lemma \ref{equal}, we get
$$
|t^2E_{\alpha, 3}(-t^\alpha)|<P^{\alpha,3-\alpha}1=t^2\frac{\Gamma(1)}{\Gamma(3)}=\frac{t^2}{2}.
$$
Therefore
$$
|E_{\alpha, 3}(-t^\alpha)|<\frac{1}{2}, \quad t>0.
$$

Now we show that $E_{\alpha, 3}(-t^\alpha)>0$. Since $\beta=3>f(\alpha)$, then from Lemma \ref{theor} we have $E_{\alpha, 3}(-t^\alpha)\neq 0$, and therefore $E_{\alpha, 3}(-t^\alpha)$ function keeps its sign for all  $t\geq 0$. But we know that  $E_{\alpha, 3}(0 )=\frac{1}{2}>0$ and $E_{\alpha, 3}(-t^\alpha)\in C[0, \infty)$. Hence $E_{\alpha, 3}(-t^\alpha)>0$ for all $t\geq 0$. Lemma \ref{lem3} is proved.

\begin{lem}\label{lem4}
Let $\alpha\in(1,2)$. Then there exists a number $C_0>0$ such that for all $t>0$ the following estimate holds:
$$
(E_{\alpha,2}(-t^\alpha))^2+E_{\alpha,3}(-t^\alpha)(1-E_{\alpha,1}(-t^\alpha))>C_0.
$$
\end{lem}

\emph{Proof.} We have:
$$
(E_{\alpha,2}(-t^\alpha))^2+E_{\alpha,3}(-t^\alpha)(1-E_{\alpha,1}(-t^\alpha))\geq E_{\alpha,3}(-t^\alpha)(1-E_{\alpha,1}(-t^\alpha)).
$$

According to Lemma \ref{lem1} and Lemma \ref{lem3}, there exists a number $C_0$, such that
$$
E_{\alpha,3}(-t^\alpha)(1-E_{\alpha,1}(-t^\alpha))>C_0.
$$
Lemma \ref{lem4} is proved.

\section{Formal solution of the problem (\ref{prob1})-(\ref{ConIn})}

Let the non-local problem (\ref{prob1})-(\ref{ConIn}) has a unique solution. Then since the system $\{v_k\}$ is complete in $H$, the solution has the form:
\begin{equation}\label{eq1}
    u(t)=\sum\limits_{k=1}^\infty T_k(t)v_k.
\end{equation}
If we multiply both sides of this equality scalarly by $v_j$, then from the orthonormality of the system of eigenfunctions \{$v_k$\}
we obtain the equalities $T_j(t)=(u(t),v_j)$. 

 Substituting (\ref{eq1}) into equation (\ref{prob1}) and using the conditions (\ref{ConB}) and (\ref{ConIn}), we have the following problem:
\begin{equation}\label{yorm}
    D_t^\alpha T_k(t)+\nu_k^2 T_k(t)=0,
\end{equation}

\begin{equation}\label{con1}
    T_k(0)=T_k(T),
\end{equation}
and
\begin{equation}\label{con2}
    \int\limits_0^T T_k(t)dt=\varphi_k.
\end{equation}
The solution to the equation (\ref{yorm}) has the form (see, for example \cite{Klb}, p. 231.)

\begin{equation}\label{sol1}
    T_k(t)=a_kE_{\alpha,1}(-\nu_k^2t^\alpha)+b_ktE_{\alpha,2}(-\nu_k^2t^\alpha).
\end{equation}
To find the unknown coefficients $a_k$ and $b_k$, we use the non-local conditions (\ref{con1}) and (\ref{con2}).

Apply conditions (\ref{con1}) and (\ref{con2}) to (\ref{sol1}), to get:
\begin{equation}\label{eq3}
\left\{
\begin{aligned}
&a_k=a_kE_{\alpha,1}(-\nu_k^2T^\alpha)+b_kTE_{\alpha,2}(-\nu_k^2T^\alpha),\\
&\int\limits_0^T (a_kE_{\alpha,1}(-\nu_k^2t^\alpha)+b_ktE_{\alpha,2}(-\nu_k^2t^\alpha))dt=\varphi_k.
\end{aligned}
\right.
\end{equation}
Solving this system of equations, we will have \begin{equation}\label{koef1}
a_k=\frac{\varphi_kE_{\alpha,2}(-\nu_k^2T^\alpha)}{T((E_{\alpha,2}(-\nu^2_kT^\alpha))^2+E_{\alpha,3}(-\nu_k^2T^\alpha)(1-E_{\alpha,1}(-\nu_k^2T^\alpha)))},    
\end{equation}

\begin{equation}\label{koef2}
b_k=\frac{\varphi_k(1-E_{\alpha,1}(-\nu_k^2T^\alpha))}{T^2((E_{\alpha,2}(-\nu^2_kT^\alpha))^2+E_{\alpha,3}(-\nu_k^2T^\alpha)(1-E_{\alpha,1}(-\nu_k^2T^\alpha)))}.    
\end{equation}

Using the equalities (\ref{eq1}), (\ref{sol1}), (\ref{koef1}) and (\ref{koef2}) we get the formal solution (\ref{yechim}) for the problem (\ref{prob1}) - (\ref{ConIn}). It remains to prove that the constructed formal solution satisfies all the requirements of Definition \ref{def}, i.e. is indeed a solution to problem (\ref{prob1}) - (\ref{ConIn}). We will do this in the next section.

On the other hand, the uniqueness of the solution follows from the already established equalities (\ref{koef1}) and (\ref{koef2}). Indeed, let us show that the solution to the homogeneous problem (\ref{prob1}) - (\ref{ConIn}) with function $\varphi=0$ is identically zero. From equalities (\ref{koef1}) and (\ref{koef2}) it follows that $a_k=b_k=0$, and then all coefficients $T_k(t)$ of series (\ref{eq1}) are equal to zero. Due to the completeness of system $\{v_k\}$, it follows that $u(t)\equiv 0$.

\section{Proof of the Theorem \ref{w2}}

Let $S_j(t)$ be the partial sums of (\ref{yechim}). 
Then
$$
AS_j(t)=\sum\limits_{k=1}^j \lambda_k(a_kE_{\alpha,1}(-\nu_k^2t^\alpha)+b_ktE_{\alpha,2}(-\nu_k^2t^\alpha))v_k.
$$
By Parseval equality:
$$
||AS_j(t)||^2=\sum\limits_{k=1}^j \lambda_k^2 |a_kE_{\alpha,1}(-\nu_k^2t^\alpha)+b_ktE_{\alpha,2}(-\nu_k^2t^\alpha)|^2\leq
$$
\begin{equation}\label{baho1}
\leq C\sum\limits_{k=1}^j \lambda_k^2 |a_kE_{\alpha,1}(-\nu_k^2t^\alpha)|^2+C\sum\limits_{k=1}^j \lambda_k^2|b_ktE_{\alpha,2}(-\nu_k^2t^\alpha)|^2.   
\end{equation}
Let us estimate the following two terms separately:
$$
I_1=|a_kE_{\alpha,1}(-\nu_k^2t^\alpha)|=\left|\frac{\varphi_kE_{\alpha,2}(-\nu_k^2T^\alpha)E_{\alpha,1}(-\nu_k^2t^\alpha)}{T((E_{\alpha,2}(-\nu^2_kT^\alpha))^2+E_{\alpha,3}(-\nu_k^2T^\alpha)(1-E_{\alpha,1}(-\nu_k^2T^\alpha)))} \right|
$$
and
$$
I_2=|b_ktE_{\alpha,2}(-\nu_k^2t^\alpha)|=\left|\frac{\varphi_kt(1-E_{\alpha,1}(-\nu_k^2T^\alpha))E_{\alpha,2}(-\nu_k^2t^\alpha)}{T^2((E_{\alpha,2}(-\nu^2_kT^\alpha))^2+E_{\alpha,3}(-\nu_k^2T^\alpha)(1-E_{\alpha,1}(-\nu_k^2T^\alpha)))}\right|.
$$

To estimate $I_1$, we apply Lemma \ref{lem1}, Lemma \ref{lem2} and Lemma \ref{lem4}. Then

\begin{equation}\label{es1}
I_1\leq\frac{|\varphi_k|}{T}\frac{1}{C_0}\leq CT^{-1}|\varphi_k|.
\end{equation}
Similarly
\begin{equation}\label{es2}
I_2\leq \frac{t|\varphi_k|}{T^2}\frac{1}{C_0}\leq CT^{-2}t|\varphi_k|.
\end{equation}
Using estimates (\ref{es1}) and (\ref{es2}) we obtain:
\begin{equation}\label{baho2}
    ||AS_j(t)||^2\leq C^2T^{-2}\sum\limits_{k=1}^j \lambda_k^2 |\varphi_k|^2+C^2T^{-4}t^2\sum\limits_{k=1}^j \lambda_k^2|\varphi_k|^2.
\end{equation}
Therefore, if $\varphi\in D(A)$ then
$$
C^2T^{-2}\sum\limits_{k=1}^j \lambda_k^2 |\varphi_k|^2+C^2T^{-4}t^2\sum\limits_{k=1}^j \lambda_k^2|\varphi_k|^2 \leq const.
$$
Thus $Au(t)\in C([0,T];D(A))$.

Now we show that the following sums $D_t^\alpha S_j(t)$ converge uniformly in $t\in (0,T)$. To do this, first consider the sums
$$
(I+A)^{-1} AS_j(t)=\sum\limits_{k=1}^j \frac{\lambda_k}{1+\lambda_k}(a_kE_{\alpha,1}(-\nu_k^2t^\alpha)+b_ktE_{\alpha,2}(-\nu_k^2t^\alpha))v_k.
$$
By Parseval equality:
$$
||(I+A)^{-1} AS_j(t)||^2=\sum\limits_{k=1}^j \frac{\lambda_k^2}{(1+\lambda_k)^2} |a_kE_{\alpha,1}(-\nu_k^2t^\alpha)+b_ktE_{\alpha,2}(-\nu_k^2t^\alpha)|^2\leq
$$
\begin{equation}\label{baho1}
\leq C\sum\limits_{k=1}^j \frac{\lambda_k^2}{(1+\lambda_k)^2} |a_kE_{\alpha,1}(-\nu_k^2t^\alpha)|^2+C\sum\limits_{k=1}^j \frac{\lambda_k^2}{(1+\lambda_k)^2}|b_ktE_{\alpha,2}(-\nu_k^2t^\alpha)|^2.   
\end{equation}
By estimates (\ref{es1}), (\ref{es2}) and $\frac{\lambda_k}{1+\lambda_k}\leq 1$ we get:
\begin{equation}\label{baho2}
    ||(I+A)^{-1}AS_j(t)||^2\leq C^2 T^{-2}\sum\limits_{k=1}^j |\varphi_k|^2+C^2T^{-4}t^2\sum\limits_{k=1}^j |\varphi_k|^2.
\end{equation}
From this, since $\varphi\in H$ we get:
$$
C^2 T^{-2}\sum\limits_{k=1}^j  |\varphi_k|^2+C^2 T^{-4}t^2\sum\limits_{k=1}^j |\varphi_k|^2 \leq const,
$$
Therefore $ (I+A)^{-1} Au(t)\in C((0,T);H)$. Now applying the obvious equality $D_t^\alpha u(t)=-\nu^2 (I+A)^{-1} Au(t)$, which follows from the commutativity of the corresponding operators, we obtain $D_t^\alpha u(t)\in C((0,T),H)$.

It remains to prove the continuity of $AD_t^\alpha u(t)$. From equality $AD_t^\alpha u(t)=-D_t^\alpha u(t) -\nu^2Au(t)$ and continuity of $D_t^\alpha u(t)$ and $Au(t)$, it follows   $AD_t^\alpha u(t)\in C((0,T),D(A))$.
Theorem is proved.

\par\bigskip\noindent
{\bf Acknowledgment.} 

The authors are grateful to Sh. A. Alimov discussions of these results.

The authors acknowledge financial support from the Innovative Development Agency under the Ministry of Higher Education, Science and Innovation of the Republic of Uzbekistan, Grant No F-FA-2021-424.

\bibliographystyle{amsplain}

\end{document}